\documentclass[11pt ]{article}

\usepackage{latexsym,amssymb}
\usepackage[applemac]{inputenc}

\newif\ifpdf
\ifx\pdfoutput\undefined
\pdffalse 
\else
\pdfoutput=1 
\pdftrue 
\fi
\ifpdf
\usepackage[pdftex]{graphicx}
\else
\usepackage[dvips]{graphicx}
\fi
 \newcommand{\Section}{\setcounter{equation}{0} \section}

\newcommand{\R} {\mathbb{R}}     
\newcommand{\eps} {\varepsilon}          %
  \newcommand{\preuve}{\noindent\textit{ Proof -~}}
 \newcommand{\findemo}{\hfill $\Box$}

\newtheorem{theo}{Theorem}[section]

\newtheorem{defi}{Definition}[section]
\newtheorem{lemme}{Lemma}[section]

\begin{document}
\ifpdf
\DeclareGraphicsExtensions{.pdf,.jpg}
\else
\DeclareGraphicsExtensions{.eps,.jpg}
\fi
\thispagestyle{empty}
   \begin{center}
{\Large{\textbf{ A new class of smoothing  methods for  mathematical programs with equilibrium constraints}}}

 M. Haddou \\
MAPMO-UMR 6628\\
Universit\'e d'orl\'eans - BP 6759 \\
45067 Orl\'eans cedex 2 \\
mounir.haddou@univ-orleans.fr

\today
\end{center}

\begin{abstract}
A  class of smoothing methods is proposed for solving mathematical programs with equimibrium constraints. We introduce new and very simple regularizations of the complementarity constraints.
Some estimate distance to optimal solution and expansions of the optimal value function are 
presented. Limited numerical experiments using SNOPT algorithm are presented to verify the efficiency of our approach.
\end{abstract}

 \Section{Introduction}
Mathematical programs with equilibrium constraints (MPECs) constitute an important class of optimization problems and pose special theoretical and numerical challenges. \\
MPECs are constrained optimization problems in which the essential constraints are defined by some parametric variational inequalities or a parametric complementarity system. 
 MPECs can be  closely related to the well-known Stackelberg game and to general bilevel programming. 
 As a result, MPECs play a very important role in many fields such as engineering
design, economic equilibrium, multilevel game, and mathematical programming
theory itself, and it has been receiving much attention in the optimization world.\\
However, MPECs are very difficult to deal with because, the feasible region and optimal solution set are  almost non convex non concave and not even connected.
Moreover, the constraints can not satisfy any standard constraint qualification such as the linear independence constraint qualification or the
Mangasarian-Fromovitz constraint qualification at any feasible point \cite{Chen,Ye}. \\
In this paper, we consider  MPECs in their standard complementarity constrained  optimization problems formulation
 \begin{equation}
 \left\{\begin{array}{lllll}\label{P}
 \min & f(x,y)\\
 {\rm s.t.} &  x\in {\cal X},y\in \R^m  , z\in\R^l, \lambda\in  \R^l,     \\
 & F(x,y)-\nabla_y g(x,y)^T\lambda=0\\
 &g(x,y)=z\\
 &z\ge 0,\lambda\ge 0, \lambda^Tz=0\\
 \end{array}\right .
 \end{equation}
 where the functions $f:\R^{n+m}\to\R$, $F:\R^{n+m}\to\R^m$ and $g:\R^{n+m}\to\R^l$ are all twice continuously differentiable and $\cal X$ is a nonempty and compact subset of  $\R^n$.\\
 {\bf Remark.} The constraints of (\ref{P}) correspond to the KKT conditions of the parametrized variational inequality
  \begin{equation}
\label{VI}
y\in C(x) {\rm \ and  \ }  (v-y)^TF(x,y)\ge 0 {\rm \ for \ all \ } v\in C(x),
 \end{equation}
 where $C(x):=\{ y\in \R^m / g(x,y)\ge 0 \}$.\\
 The negative properties of MPECs make these problems very  difficult and exclude any direct use of 
 standard non linear programming (NLP) algorithms. \\
 In this paper we propose some smoothing techniques to regularize the complementarity constraints
 and construct  relaxed  problems that are suitable for NLP algorithms.\\
 Many regularization and relaxation techniques have already been proposed, here is an incomplete list of such methods\\
\begin{equation}
 \begin{array}{lllll}
 (Reg(t)\cite{RW,Shol}) &\lambda^Tz=0 {\rm\ is \ relaxed \ to }& \lambda_i z_i\le t \quad\forall i\\
 (Regeq(t)\cite{RW,Shol})&\lambda^Tz=0 {\rm\ is \ replaced \ by }&  \lambda_i z_i= t  \quad\forall i\\
 (RegCp(t)\cite{RW,Shol})& \lambda^Tz=0{\rm\ is \ relaxed \ to}& \lambda^Tz\le t\\
 (Facc.\cite{Facc})&\lambda^Tz=0{\rm\ is \ replaced \ by}& \sqrt{( \lambda_i - z_i)^2 +4t^2}-( \lambda_i + z_i)=0\ \forall i\\
 (Entro.\cite{Bir,Fang})&\lambda^Tz=0{\rm\ is \ replaced \ by}& \displaystyle t\ln\{ e^{\frac{- \lambda_i}{t}}+ e^{\frac{ -z_i}{t}}\} =0\quad\forall i.\\
 \end{array}
 \end{equation}
 In almost all these techniques,  the constraints $\lambda_i z_i=0$ or  $\min(\lambda_i, z_i)=0$ are replaced by some smooth approximations.\\
 In our approach, the complementarity constraint is interpreted componnent-wise as:\\
 \centerline{ $\forall i,$ \quad At most one of $z_i$ or $\lambda_i$ is nonzero.}\\
 So, we  construct some parameterized real functions that satisfy:\\
 \[ (\theta_r(x)\simeq1 {\rm \ if \ } x\ne 0) {\rm \ and  \ } (\theta_r(x)\simeq 0 {\rm \ if\  } x= 0)\]
 to count nonzeros and then replace the constraint
 $$\lambda_i z_i=0$$
  by 
  $$\theta_r(\lambda_i)+\theta_r(z_i)\le 1.$$
 In section 2, we present some preliminaries and assumptions on the problem (\ref{P}) (essentially the same as in \cite{RW}) . In Section 3,  the smoothing functions and techniques are presented and many approximation and  regularity properties are proved. Section 4 is devoted to the analysis of the regularization process. The last section presents some numerical experiments concerning 
 two  smoothing functions. 
 
 \Section{Assumptions  and preliminaries} 
We essentially need the same assumptions and background as in \cite{RW}. A complete presentation 
of this background needs about 6 to 7 pages.
We will  only present in this section some definitions, known optimality conditions and constraint qualifications. For some others we will only refer readers to \cite{RW}. 
These notions will be useful in the next section. 
 \\
 The first definition concern a first order optimality condition: the strong stationarity
 \begin{defi} A feasible point $(x^*,y^*,z^*,\lambda^*)$ is strongly stationary for (\ref{P}) if $d=0$ solves 
 \begin{equation}
 \left\{
 \begin{array}{lllll}\label{ststat}
 \min & \nabla f(x^*,y^*)^Td_{x,y}\\
 {\rm s.t.} &  d_x\in Add({\cal X}(x^*)) , d_z\in\R^l, d_\lambda\in  \R^l,     \\
 &\nabla F(x^*,y^*)^Td_{x,y} -\nabla_y g(x^*,y^*)^Td_\lambda - \nabla (\nabla_y g(x^*,y^*))^Td_{x,y} =0\\
 &\nabla g(x^*,y^*)^Td_{x,y}-d_z=0\\
 &(d_z)_i=0, i\in I_z\backslash I_\lambda\\
 &(d_\lambda)_i=0, i\in I_\lambda\backslash I_z\\
 &(d_z)_i\ge 0,(d_\lambda)_i\ge 0, i\in I_z\cap I_\lambda\\
 \end{array}\right .
 \end{equation} 
 \end{defi}
 where $d=(d_x,d_y,d_z,d_\lambda)^T\in  \R^{n+m+2l} $, $I_z$ and $I_\lambda$ are the active sets at $(x^*,y^*,z^*,\lambda^*)$
 $$I_z:=\{i=1,\dots, l\vert z^*_i=0\}{\rm\  and\  }I_\lambda:=\{i=1,\dots, l\vert \lambda^*_i=0\}$$
 and $Add({\cal X}(x^*))$ is the admissible directions set defined by  
 $$Add({\cal X}(x^*)):=\{d_x\in \R^m\vert \exists r_0>0\quad  \forall 0\le r\le r_0 \quad x^*+rd_x\in {\cal X}.\}$$
{\bf Remark.} There is an other kind of stationarity (the B-stationarity) which is {\bf less} restrictive but {\bf very difficult to chek}. We prefer to not present it in this paper. These two stationarity properties are equivalent 
when the MPEC-LICQ  (defined next) is satisfied
 \begin{defi} The MPEC-LICQ  is satisfied at the point $(x^*,y^*,z^*,\lambda^*)$ if the linear independance constraint qualification (LICQ) is satisfied for the following RNLP problem at $(x^*,y^*,z^*,\lambda^*)$.
 \end{defi}
 \begin{equation}
 \left\{
 \begin{array}{lllll}\label{RNLP}
 \min & f(x,y)\\
 {\rm s.t.} &  x\in {\cal X} , z\in\R^l, \lambda\in  \R^l,     \\
 & F(x,y)-\nabla_y g(x,y)^T\lambda=0\\
 &g(x,y)=z\\
 &z_i=0, i\in I_z\backslash I_\lambda\\
 &\lambda_i=0, i\in I_\lambda\backslash I_z\\
 &z_i\ge 0,\lambda_i\ge 0, i\in I_z\cap I_\lambda\\
 \end{array}\right .
 \end{equation}  
 An other important and usefull constraint qualification is the following Mangasarian-Fromovitz one
 \begin{defi} The MPEC-MFCQ  is satisfied at the point $(x^*,y^*,z^*,\lambda^*)$ if the  Mangasarian-Fromovitz constraint qualification (MFCQ) is satisfied for the RNLP problem at $(x^*,y^*,z^*,\lambda^*)$.
 \end{defi} 
 
We will also use some Second-Order sufficient conditions namely: 
the (MPEC-SOSC) and the (RNLP-SOSC).
These two conditions (among others) are defined in \cite{RW}.  
 
\Section{The smoothing technique}
For $r>0$, we consider real functions $\theta_r: \R_+\rightarrow [0,1]$ satisfying
 \begin{equation}\label{condtheta}
\begin{array}{llll}
 (i)& \theta_r {\rm\  is\  nondecreasing,\ strictly \ concave\  and\  continuously\  differentiable,}\\
(ii)&\forall r>0, \quad \theta_r(0)=0,\\
(iii)&\forall x>0,  \quad\displaystyle \lim_{r\to 0}\theta_r(x)=1,\  {and}\\
(iv)&\displaystyle \lim_{r\to 0}\theta'_r(0)>0. 
\end{array}
 \end{equation}
We will present some interesting examples of such functions after the following approximation result
 \begin{lemme}\label{lem1} For any $\eps >0$, and $ x,y\ge 0$, there exists $ r_0>0$ such that
$$ \forall r\le r_0, \qquad  (\min (x,y)=0) \Longrightarrow\  (\theta_r(x)+\theta_r(y)\le 1)\Longrightarrow\  (\min (x,y)\le \varepsilon). $$
\end{lemme}
\preuve 
The first property  is obvious since $\theta_r(0)=0$ and $\theta_r\le 1$.\\
Using assumption $(iii)$ for $x=\varepsilon$, we have
\[
\forall \alpha>0, \quad \exists r_0>0/ \quad \forall r\le r_0 \qquad 1-\theta_r(\varepsilon)<\alpha,
\]
so that, if we suppose that $\min (x,y)>\varepsilon$, assumption $(i)$ gives
\[
\theta_r(x)+\theta_r(y)>2\theta_r(\varepsilon)> 2(1-\alpha).
\]
Then if we choose  $\alpha<\frac{1}{2}$, we obtain that $\theta_r(x)+\theta_r(y)>1$.
\findemo\\
This first approximation result can be improved for some interesting choices of the smoothing functions $\theta_r$
\[
 \begin{array}{lllll}
(\theta^1_\cdot) \qquad&\displaystyle \theta^1_r(x)=\frac{x}{x+r}\\
 (\theta^{W_k}_\cdot)\qquad &\displaystyle \theta^{W_k}_r(x)=1-e^{-(\frac{x}{r})^k}\quad {\rm \ for\ } k>0\\
 (\theta^{\log}_\cdot)\qquad&\displaystyle \theta^{\log}_r(x)=\frac{\log(1+x)}{\log(1+x+r)}\\

 \end{array}
 \]
 We will also consider the general class ${\Theta}^{\ge 1}$ of functions 
\[
 \begin{array}{lllll} 
 (\theta^{\ge 1}_\cdot)\qquad&  {\rm verifying\ (i-iv)\ and\ } \theta^{\ge 1}_\cdot\ge \theta^1_\cdot\\
 \end{array}
 \] 
 {\bf  Remark.} The functions $ \theta^{W_k}_\cdot$ are the density functions of Weibull distributions, when $k=1$, the obtained smoothing method corresponds (with slight modifications) to the inequality entropic regularization \cite{Bir}. Simple comparison calculus proove that $ \theta^{\log}_\cdot$ and $ \theta^{W_k}_\cdot$ for ($0<k\le 1$) belong to the class of functions $\Theta^{\ge 1}$.
 \begin{lemme}\label{lem2} we have
 \[
 \begin{array}{lllll}
(i)&\forall x\ge 0,\ \forall y\ge 0 \qquad&\theta^1_r(x)+\theta^1_r(y)\le 1 \Longleftrightarrow x\cdot y\le r^2, and\\
(ii)&\forall x\ge 0,\ \forall y\ge 0 \qquad &x\cdot y=0\Longrightarrow\theta^{\ge 1}_r(x)+\theta^{\ge 1}_r(y)\le 1 \Longrightarrow x\cdot y\le r^2.\\

 \end{array}
 \]
\end{lemme}
\preuve 
(i) We have 
\[
\displaystyle \theta^1_r(x)+\theta^1_r(y)=\frac{2xy+rx+ry}{xy+rx+ry+r^2},
\]
so that
\[
 \begin{array}{lllll}
  \theta^1_r(x)+\theta^1_r(y)\le 1 &\Longleftrightarrow &2xy+rx+ry\le xy+rx+ry+r^2\\
 &\Longleftrightarrow &x\cdot y\le r^2.\\

 \end{array}
\]
The first part of (ii) follows obviously from  Lemma \ref{lem1} and the second one is a direct consequence of (i) since
\[
\theta^{\ge 1}_r(x)+\theta^{\ge 1}_r(y)\le 1 \Longrightarrow\theta^1_r(x)+\theta^1_r(y)\le 1.
\]
\findemo \\
Using any function $\theta_r$ satisfying (\ref{condtheta}), we obtain the relaxed following 
problem for (\ref{P})
 \begin{equation}
 \left\{
 \begin{array}{lllll}\label{Pr}
 \min & f(x,y)\\
 {\rm s.t.} &  (x,y)\in {\cal X} , z\in\R^l_+, \lambda\in  \R^l_+,   e \in  \R^l_+ \\
 & F(x,y)-\nabla_y g(x,y)^T\lambda=0\\
 &g(x,y)=z\\
 & \theta_r(\lambda_i)+\theta_r(z_i)+ e_i= 1, \quad \forall i\in\{1,...,l\}.\\
 \end{array}\right.
 \end{equation}
{\bf Remarks.}  (i) By choosing some particular smoothing functions (ex. $\theta^{W_k}_r$), the nonnegativity constraints on $\lambda$ and $z$ become implicite and can be removed from the definition of (\ref{Pr}).(This can have an important impact in practice.)\\
(ii) Under some classical assumptions, as in \cite{Facc} we can easily prove that  the jacobian of equality constraints (with respect to $(y,z,\lambda)$ is nonsingular. This property is useful in practice since standard NLP algorithms use Newton-type to solve systems of nonlinear equations corresponding to this jacobian.
  \begin{lemme}\label{lem3} If $g$ is concave with respect to $y$ and $F$ is uniformly strongly monotone with respect to $y$, then for every nonnegative $r$ and every feasible point $(x,y,z,\lambda,e)$ of problem (\ref{Pr}), 
  the jacobian of equality constraints (with respect to $(y,z,\lambda)$) is nonsingular.
  \end{lemme}
 \preuve  
  Using the assumptions \ref{condtheta} (i) and (iv), the proof is exactely the same as in \cite{Bir} or \cite{Facc}.
  \findemo 
  
  Problem  (\ref{Pr}) may be viewed as a perturbation of  (\ref{P}). Previous lemmas prove that  (\ref{Pr}) is in fact some tight relaxation of  (\ref{P}). However this perturbation is not continuous on the parameter $r$ so that any direct use of perturbation results such that \cite{BS} is impossible. \\
Fortunately,  Lemma \ref{lem2} proves that for the particular smoothing function $\theta^1_\cdot $, the corresponding relaxed problem (\ref{Pr}) is equivalent to $(Reg({\bf r^2})) $ in \cite{RW}. We can then benefit from the theoretical results in \cite{RW}. \\
 The following results provide, in the case of the {\bf $\theta^1_\cdot $} function, some distance estimate  between  
  solution of  (\ref{Pr}) and solution of (\ref{P}). These results correspond to applications of [\cite{BS}, Theorem 5.57, Theorem 4.55 and Lemma 4.57] and can be found with complete proofs in \cite{RW}. We just state them in our context and  add the optimal value expansion. 
  \begin{theo}\label{theo1} Suppose that $X^*=(x^*,y^*,z^*,\lambda^*)$ is a strongly stationary point of  (\ref{P}) at which MPEC-MFCQ and MPEC-SOSC are satisfied. Then there are positive constants $\alpha$, $\bar r$, and $M$ such that for all $r\in (0,\bar r]$, the global solution $X(r)$ of the localized problem (\ref{Pr}) with the additional ball constraint $\Vert X-X^*\Vert \le \alpha$ that lies closest to $X^*$ satisfies $\Vert X(r)-X^*\Vert \le M.r$. Furthermore the optimal value  $v_r$of  (\ref{Pr}) has an expansion of the form
  $$v_r= v^0+\frac{1}{2}.a.r^2 +o(r^2)$$
  where $v^0$ is the optimal value of (\ref{P}) and  $a$ is the optimal value of an auxiliary quadratic problem\cite{BS}.

 \end{theo}
   \begin{theo}\label{theo2} Suppose that $X^*=(x^*,y^*,z^*,\lambda^*)$ is a strongly stationary point of  (\ref{P}) at which MPEC-LICQ and RNLP-SOSC are satisfied. Then there are positive constants $\alpha$, $\bar r$, and $M$ such that for all $r\in (0,\bar r]$, the global solution $X(r)$ of the localized problem (\ref{Pr}) with the additional ball constraint $\Vert X-X^*\Vert \le \alpha$ that lies closest to $X^*$ satisfies $\Vert X(r)-X^*\Vert \le M.r^2$. Furthermore the optimal value $v_r$ of  (\ref{Pr}) has an expansion of the form
  $$v_r\le v^0+b.r^2 +O(r^4)$$
  where $v^0$ is the optimal value of (\ref{P}) and  $b$ is the optimal value of an auxiliary linearized problem\cite{BS}.
 \end{theo}
 For functions  of the general class $\Theta^{\ge 1}$, the corresponding feasible sets satisfy
 \[
{\cal F}_P \subset {\cal F}_{\theta^{\ge 1}_\cdot} \subset {\cal F}_{\theta^1_\cdot} 
 \]
 where ${\cal F}_P$, ${\cal F}_{\theta^{\ge 1}_\cdot} $ and ${\cal F}_{\theta^1_\cdot}$ are respectively the feasible set of problem  (\ref{P}) and (\ref{Pr}) for the corresponding $\theta_r$ function.\\
These inclusions prove  that the optimal value expansions given in  Theorem \ref{theo1} and Theorem \ref{theo2} are still valid under the same assumptions.
   \begin{theo}\label{corro} When using functions  $\theta^{\ge 1}_\cdot$, under the same assumptions of Theorem\ref{theo1} (resp. Theorem\ref{theo2})
  the optimal value $v_r$ of  (\ref{Pr}) has an expansion of the form
  $$ v_r\le  v^0+\frac{1}{2}.a.r^2 +o(r^2) \quad ({\rm resp.\ }\quad v_r\le v^0+b.r^2 +O(r^4))$$
 \end{theo}


 \Section{Numerical resuts}
For two different smoothing functions, we present some numerical results using the SNOPT \cite{Gill}
nonlinear programming algorithm on the AMPL \cite{AMPL} optimization plateform. Our aim is just to verify the qualitative numerical efficiency of our approach. We consider a  subset of the MACMPEC \cite{MACMPEC} test problems with known optimal values and solutions (a large part of these test problems were used by \cite{Bir, Facc} in their numerical experiments) .\\ 
We choose the two functions
$$\theta^1_r(x)=\frac{x}{x+r}$$
and 
$$\displaystyle\theta^{W_1}_r(x)=1-e^{-\frac{x}{r}}.$$
The first function  has (in our analysis) the best theoretical results and corresponds ¨in some way¨to 
the regularization studied in  \cite{Shol,RW}. While the second one corresponds to the enropic regularization \cite{Bir,Fang}.\\
In our experiments, we made a logarithmic scaling for these two functions to bound their gradients.
Each constraint 
$$\theta_r(\lambda_i)+\theta_r(z_i)+ e_i= 1$$
is  in fact replaced by the following inequality
$$\displaystyle r^2\ln \left( \frac{r}{\lambda_i+r}+ \frac{r}{z_i+r}\right) \ge 0,$$
in the case of the $\theta^1_r$ function and 
$$\displaystyle r\ln \left( e^{-\frac{\lambda_i}{r}}+ e^{-\frac{z_i}{r}}\right) \ge 0.$$
in the case of the $\theta^{W_1}_r$ function.\\
The two following tables give for each considered problem and for different starting points, the used value of the parameter $r$, the optimal value and solution obtained when using each of the two smoothing functions.
The tables report also different informations concerning the computational effort of the solver SNOPT. $itM$ and $itm$ correspond to the total number of major and minor iterations numbers \cite{Gill}.  The total number of objective function evaluations is given in $(Obj.)$. $(grad.)$ corresponds to the total number of objective function gradient evaluations. $(constr.)$ and $(jac.)$ give respectively the total number of constraints and  constraints gradient evaluations. \\

\begin{center}
{\scriptsize
\begin{tabular}{llllllllllll}    \hline
Problem &r& Start & Obj.val. & Opt.x& (itM,itm)&Obj. & grad& constr.& Jac \\ \hline
Bard1 &1.e-2 &no& 17 &(1,0)&(5,8)& 9& 8 & 9&8\\ \hline
Df1 &1.e-3 &no& 0 &(1,0) &(1,1)& 3 & 2 & 3&2\\ \hline
Gauvin &1.e-2  &no& 20 &(2,14) &(4,11)& 7 & 6& 7&6\\ \hline
jr1& 1.e-2&no& 0.5 &(0.5,0.5) &(6,3)& 9 & 8 & 9 &8\\ \hline

Gnash10 &1.e-5 &gnash10.dat& -230.8232 &47.036&(17,46)& 21 & 20 & 21 &20\\ \hline
Gnash11 &1.e-4 &gnash11.dat& -129.9119 & 34.9942&(20,50)& 18 & 17 & 21 &20\\ \hline
Gnash12 &1.e-4 &gnash12.dat& -36.93311 &18.1332&(24,51)& 27 & 26 & 27 &26\\ \hline
Gnash13 &1.e-2 &gnash13.dat& -7.061783&7.55197&(14,56)& 20 & 19 & 23 &22\\ \hline
Gnash14 &1.e-3 &gnash14.dat& -0.179046 &1.06632&(14,46)& 18 & 17 & 21 &20\\ \hline

Scholtes1&1.e-1 &1& 2&0 &(9,10)& 14 & 13 & 14 &13\\ \hline
Bilevel1 &1.e-2 &(25,25)& 5&(25,30)&(3,11)& 0& 0& 9 &8\\ \hline
 & &(50,50)& 5&(25,30)  &(0,6)& 0 & 0 & 2&1 \\ \hline

Nash1 & 1.e-1&(0,0)& 1.61e-14&(9.996,4.999) &(13,42)& 25 & 24 & 25&24 \\ \hline
 & &(5,5)& 1.60e-18 &(9.313,5.686) &(10,33)& 32 & 31 & 32&31 \\ \hline
 & &(10,10)& 1.46e-14&(9.092,5.901) &(16,38)& 34& 33 & 34&33 \\ \hline
 & &(10,0)& 3.56e-24 &(9.999,4.999) &(12,34)& 28 & 27 & 28 &27\\ \hline
 & &(0,10)& 9.03e-22 &(9.999,4.999)&(14,41)& 31 & 30 & 31 &30\\ \hline

Bilevel2 & 1.e-4&(0,0,0,0)&   -6600&(6.441,4.863,12.559,16.137) &(6,43)& 9& 8 & 9&8 \\ \hline
 & &(0,5,0,20)& -6600&(6.575,5,12.425,16)&(6,50)& 9 & 8 & 9&8 \\ \hline
 & &(5,0,15,10)& -6600&(6.837,12.162,16)&(5,36)& 7& 6 & 7&6 \\ \hline
 & &(5,5,15,15)& -6600&( 4.892,3.373,14.107,17.627)&(5,35)& 7 & 6 & 7&6 \\ \hline
 & &(10,5,15,10)& -6600&(8.014,4.971,10.986,16.029)&(5,38)& 7 & 6 & 7&6 \\ \hline
Bilevel3 &1.e-4 &(0,0) & -12.6787 &(0,2)&(9,23)& 12 & 11 & 12&11 \\ \hline
 & &(0,2) &-12.6787 &(0,2)&(15,27)& 32 & 31 & 32 &31\\ \hline
 & &(2,0) &-10.36  &(2,0)&(01,06)& 3 & 2 & 3 &2\\ \hline

desilva &1.e-3 &(0,0)& -1 &(0.5,0.5) &(4,11)& 6& 5 & 6 &5\\ \hline
              &          &(2,2)& -1 &(0.5,0.5) &(3,9)& 5 & 4 & 5 &4\\ \hline
Stack.1 &1.e-2 &0    &  -3266.6666 &93.3333 &(4,9)& 6 & 5 & 6 &5\\ \hline
              &          &100&  -3266.6666 &93.3333 &(3,3)& 5 & 4 & 5 &4\\ \hline
              &          &200&  -3266.6666 &93.3333 &(7,5)& 19 & 18 & 19 &18\\ \hline

\end{tabular}}
\vskip 1.cm
Table1:  using the $\theta^1_r$ smoothing function
\end{center}
\begin{center}
{\scriptsize
\begin{tabular}{llllllllllll}    \hline
Problem &r& Start & Obj.val. & Opt.x& (itM,itm)&Obj. & grad& constr.& Jac \\ \hline
Bard1 &1.e-2 &no& 17 &(1,0)&(13,8)& 16 & 15 & 16&15\\ \hline
Df1 &1.e-3 &no& 0 &(1,0) &(1,2)& 3 & 2 & 3&2\\ \hline
Gauvin &1.e-2  &no& 20 &(2,14) &(5,12)& 7 & 6& 7&6\\ \hline
jr1& 1.e-2&no& 0.5 &(0.5,0.5) &(13,4)& 16 & 15 & 16 &15\\ \hline

Gnash10 &1.e-3 &gnash10.dat& -230.8232 &47.036&(17,63)& 19 & 18 & 21 &20\\ \hline
Gnash11 &1.e-3 &gnash11.dat& -129.9119 & 34.9942&(15,48)& 18 & 17 & 18 &17\\ \hline
Gnash12 &1.e-1 &gnash12.dat& -36.93311 &18.1332&(15,43)& 19 & 18 & 19 &18\\ \hline
Gnash13 &1.e-1 &gnash13.dat& -7.061783&7.55197&(23,69)& 30 & 29 & 30 &29\\ \hline
Gnash14 &1.e-3 &gnash14.dat& -0.179046 &1.06633&(22,38)& 27 & 26 & 27 &26\\ \hline

Scholtes1&1.e-1 &1& 2&0 &(11,11)& 16 & 15 & 16 &15\\ \hline

Bilevel1 &1.e-2 &(25,25)& 5&(25,30)&(3,11)& 0& 0& 9 &8\\ \hline
 & &(50,50)& 5&(25,30)  &(0,6)& 0 & 0 & 2&1 \\ \hline

Nash1 & 1.e-1&(0,0)& 7.27e-14&(9,6) &(9,16)& 12 & 11 & 12&11 \\ \hline
 & &(5,5)& 4.25e-18 &(10,5) &(6,16)& 10 & 9 & 10&9 \\ \hline
 & &(10,10)& 1.09e-11&(9,6) &(13,25)& 21& 20 & 21&20 \\ \hline
 & &(10,0)& 1.27e-13 &(9.355,5.645) &(16,34)& 27 & 26 & 27 &26\\ \hline
 & &(0,10)& 3.29e-15 &(9.396,5.604)&(7,16)& 11 & 10 & 11 &10\\ \hline

Bilevel2 & 1.e-1&(0,0,0,0)&   -6600&( 4.851,5,14.149,16) &(6,34)& 8& 7 & 8&7 \\ \hline
 & &(0,5,0,20)& -6600&(5.195,5,13.805,16)&(5,40)& 7 & 6 & 7&6 \\ \hline
 & &(5,0,15,10)& -6600&(6.099,4.834,12.901,16.166)&(5,42)& 7& 6 & 7&6 \\ \hline
 & &(5,5,15,15)& -6600&(4,1.714,15,19.286)&(5,45)& 7 & 6 & 7&6 \\ \hline
 & &(10,5,15,10)& -6600&(7.724,5,11.276,16)&(5,50)& 7 & 6 & 7&6 \\ \hline
Bilevel3 &1.e-1 &(0,0) & -12.6787 &(0,2)&(22,38)& 29 & 28 & 29&28 \\ \hline
 & &(0,2) &-12.6787 &(0,2)&(34,61)& 56 & 55 & 56 &55\\ \hline
 & &(2,0) &-10.36  &(2,0)&(01,06)& 3 & 2 & 3 &2\\ \hline

 desilva &1.e-2 &(0,0)& -1 &(0.5,0.5) &(5,9)& 8& 7 & 8 &7\\ \hline
              &          &(2,2)& -1 &(0.5,0.5) &(6,14)& 9 & 8 & 9 &8\\ \hline
 Stack.1 &1.e-2 &0    &  -3266.6666 &93.3333 &(4,4)& 6 & 5 & 6 &5\\ \hline
              &          &100&  -3266.6666 &93.3333 &(3,5)& 5 & 4 & 5 &4\\ \hline
              &          &200&  -3266.6666 &93.3333 &(11,4)& 14 & 13 & 14 &13\\ \hline
\end{tabular}}
\vskip 1.cm
Table2: using the inequality enropic approach ($\theta^{W_1}_r$)

\end{center}

\Section{Conclusion} 
We introduced a new regularization scheme for mathematical programs with complementarity constrains. Our approach  is very simple and quite different from existing techniques for the same class of problems. The obtained regularized problems are now suitable for standard NLP algorithms.
 These regularizations have different theoretical sensivity and regularity properties.
 The limited numerical experiments give very promising results (comparable to those of \cite{Bir}) and suggest to make  real investigations
 on functions of  the class $\theta^{W_k}_\cdot$.
 Therefore, we hope that some of our smoothing functions will correspond to simple and efficient
 algorithms  for the solution of real-world MPECs and Bilevel programs.

 \bibliographystyle{plain}

   \end{document}